\documentclass{amsart}
\usepackage{fullpage}
\usepackage{graphicx}
\input pictex.sty
\newtheorem{thm}{Theorem}[section]

\newtheorem{lem}[thm]{Lemma}

\theoremstyle{definition}

\theoremstyle{remark}
\newtheorem{rem}[thm]{Remark}
\numberwithin{equation}{section}

\newcommand{\To}{\longrightarrow}

\begin{document}

\title{Combinatoric Results for graphical enumeration and the higher
  Catalan numbers}%
\author{Virgil U. Pierce}%
\address{Dept. of Math, The Ohio State University}%
\email{vpierce@math.ohio-state.edu}

\thanks{This work was 
supported in part by NSF grant no. 0073087 and
no. 0412310.}%
\subjclass{}%
\keywords{}%

%\date{}%
%\dedicatory{}%
%\commby{}%
% ----------------------------------------------------------------
\begin{abstract}
We summarize some combinatoric problems solved by the higher Catalan
numbers.  These problems are generalizations of the combinatoric problems
solved by the Catalan numbers.  
The generating function of the higher Catalan numbers appeared
recently as an auxiliary function in enumerating maps and explicit
computations of the asymptotic expansion of the partition function of
random matrices in the unitary ensemble case.  
We give combinatoric proofs of the 
formulas for the number of genus 0 and genus 1 maps.
\end{abstract}
\maketitle
%\tableofcontents
% ----------------------------------------------------------------

\section{Higher Catalan Numbers}

\begin{center}
\begin{figure}[h]
\resizebox{3cm}{3cm}{\includegraphics{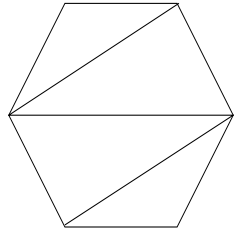}}
 \caption{\label{fig1} A polygon with $4+2 = 6$ sides divided into $4$
 triangles using edges and diagonals}
\end{figure}\end{center}

The Catalan numbers solve a number of classical combinatoric problems
such as the ``Euler Polygon Division Problem'': how many ways are
there to divide a marked polygon with $j+2$ sides into triangles using
edges and diagonals \cite{BB, Do, Go56, Ho,Po} (see figure \ref{fig1}).  
They count the number of right-left
paths along a 1-Dimensional integer lattice which stay to the right of 0 and
go from 0 to 0 in $2j$ steps; equivalently they count Dyck paths from
$(0,0)$ to $(2j, 0)$ \cite{Ba, Br, MW, Sa}.  
They count the number of ways for $2j$
customers to line up, with $j$ customers having only a 2-dollar bill and $j$
customers having only a 1-dollar bill, to purchase 1-dollar widgets, so that
each customer receives exact change.    They count the
number of non-crossing handshakes possible across a round table between $n$
people \cite{CG}.  

In this paper we will explore a generalization of the Catalan numbers,
the higher Catalan numbers.
We will show that this generalization solves enumerative problems that
are natural generalizations of the problems solved by the Catalan
numbers.  We will then highlight their appearance in recent results on
map enumeration problems.

Let
\begin{equation} \label{z(s)_taylor}
z(s) = 1 + \sum_{j=1}^\infty \zeta_j^{(\nu)} s^j
\end{equation}
be defined implicitly as the solution to the algebraic equation 
\begin{equation} \label{constraint}
s z(s)^\nu - z(s) + 1 = 0,
\end{equation}
which is regular at $s=0$.  
This generating function was presented in 
\cite{PS, YY}.
It is known  that the $\zeta_j^{(\nu)}$ are  (see \cite{Kl, YY})
\begin{equation}\label{High_Cat}
\zeta_j^{(\nu)} = \frac{1}{j} \binom{\nu j}{j-1} .
\end{equation}
The numbers (\ref{High_Cat}) are called the higher Catalan numbers
 and solve a number of combinatoric problems.

\begin{center}\begin{figure}[h]
\resizebox{3cm}{3cm}{ 
\includegraphics{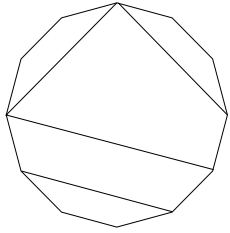} } 
\caption{\label{fig2} A polygon with $2\cdot 5 + 2 = 12$ sides divided
  in to $5$ squares, i.e. $\nu=3$, using edges and diagonals}
\end{figure}\end{center}

The $\zeta_j^{(2)}$
are the Catalan numbers.  
  The
higher Catalan numbers count the number of ways to divide a marked $((\nu-1) j
+ 2)$-sided polygon 
into $(\nu +1)$-sided polygons using edges and diagonals as in figure 
\ref{fig2}.  

They count the
number of paths along the 1-Dimensional integer lattice which move to the
right 1 unit per step, to the left $(\nu-1)$ units per step, stay to the right
of $0$ and go from 0 to 0 in $\nu j$ steps \cite{Sa}. 
They count the number of higher Dyck paths: paths which go from
$(0,0)$ to $(\nu j, 0)$ along the 2-D integer lattice with steps $(1,
1)$ and $(1, -(\nu-1))$ which stay above the $x$-axis. 
They count the number of ways for
$(\nu-1)j$ customers with 1-dollar bills and $j$ customers with $\nu$-dollar
bills to line up to buy 1-dollar widgets so that each customer receives
correct change.  
They count the number of non-crossing handshakes possible across a round table
between $n$ $\nu$-handed beings.  

\begin{rem}
More generally we call generalized Dyck paths 
walks on the 2-D integer lattice restricted to the upper half plane
with steps coming from a set of integer vectors.  The
generating function of the number of generalized Dyck paths satisfies
a system of algebraic equations \cite{LaYe90}.
\end{rem}

In addition to these combinatoric problems the 
 coefficients of $z(s)$ satisfy the recursion relation 
\begin{equation}\label{recur}
\zeta_j^{(\nu)} = \sum_{j_1 + j_2 + \dots + j_\nu = j-1} \zeta_{j_1}^{(\nu)}
\zeta_{j_2}^{(\nu)} \dots \zeta_{j_\nu}^{(\nu)}. 
\end{equation}

These numbers appeared in counting labeled maps using the
partition function of random matrices \cite{EMP}.  We found, in
solving the map enumeration
 problem for some simple cases, that the generating
function could be written as a function
of the generating function for the higher Catalan numbers.  

A map $D$ on a compact, oriented connected surface $X$ is a pair $D = ( K(D),
[\iota])$ where
\begin{enumerate}
\item $K(D)$ is a connected 1-complex;
\item $[\iota]$ is an isotopical class of inclusions $\iota:K(D) \To X$;
\item the complement of $K(D)$ in $x$ is a disjoint union of open cells
  (faces);
\item the complement of $K_0(D)$ (vertices) in $K(D)$ is a disjoint union of
  open segments (edges).
\end{enumerate}
Maps can be thought of as ribbon graphs (or fattened graphs) 
embedded on an oriented surface \cite{BIZ, Wi, Zv}.  

Define $e_g(t)$ to be the generating function of the number
of connected maps of degree $2\nu$ embedded in a genus $g$ surface, in
the following sense 
\begin{equation*}
e_g(t) = \sum_{j=1}^\infty \kappa_g^{(2\nu)}(j) \frac{ (-t)^j}{j!},
\end{equation*}
where $\kappa_g^{(2\nu)}(j)$ is the number of planar maps with $j$
vertices of degree $2\nu$ embedded in a genus $g$ surface.   

Define the class of functions: \textit{iterated integrals of rational
  functions of $z$}, or \textit{iir} of $z$, to be functions
  which are found by taking a finite number of anti-derivatives of
  a rational function of $z$.  The particular
  subclass we are going to be concerned with are those \textit{iir}'s 
  which have
  singularities only at $z =0, 1,$ or $\nu/(\nu-1)$.  One may check
  that this subclass of functions is closed under anti-differentiation.

To illustrate the results of \cite{EMP} we present the explicit result
we found for $g=0$.  
Let
\begin{equation*}
c_\nu = 2\nu \binom{2\nu-1}{\nu-1}.
\end{equation*}
\begin{thm}[Ercolani-McLaughlin-Pierce]\label{theorem_1}
\begin{equation*}
e_0(-t) = \mu_\nu (z(c_\nu t)-1) (z(c_\nu t) - r_\nu ) + \frac{1}{2}
\log(z(c_\nu t))\,,
\end{equation*}
where
\begin{equation*}
\mu_\nu = \frac{(\nu-1)^2}{4\nu (\nu+1)} 
\end{equation*}
and
\begin{equation*}
r_\nu = \frac{3(\nu+1)}{\nu-1}.
\end{equation*}
More generally we find that $e_g(-t)$ is an \textit{iir}
function of $z(c_\nu t)$ with singularities only possible at $z=0, 1$,
or $\nu/(\nu-1)$.    
\end{thm}

We will now outline the approach we used to prove this result.
Define the partition function of random matrices  
\begin{equation} \label{partfun}
Z_N(t_1, t) = \int_{\mathcal{H}_N} \exp\left( -N \mbox{Tr}\left( t_1 M +
\frac{1}{2} M^2 + t M^{2\nu} \right) \right) dM,
\end{equation}
where the integral is taken over the space of $N \times N$ Hermitian
matrices, and $dM$ is the product of Lebesgue measures on the independent
variables in $M$:
\begin{equation*}
dM = \left[ \prod_{i<j}
  d\left(\mbox{Re}(M_{ij})\right)d\left(\mbox{Im}(M_{ij})\right) \right]
  \left[ \prod_i dM_{ii} \right] \,.
\end{equation*}  

Let 
\[ L = \begin{pmatrix} b_1 & a_1 & 0 & 0 & 0 & \dots \\
  a_1 & b_2 & a_2 & 0 & 0 & \dots \\
0 & a_2 & b_3 & a_3 & 0 &\dots \\
\vdots &  & \ddots & \ddots &\ddots &  
\end{pmatrix}\,.
\]
The Toda lattice hierarchy is given by the differential equations
\[ \frac{dL}{dt_j} = \left[ B_j, L\right]\,, \]
where  $ B_j = (L^j)_+ - (L^j)_- $, and where $(\cdot)_\pm$ indicates projection
onto the upper (resp. lower) triangular parts.  
This hierarchy is an integrable hierarchy, possessing a complete family of
independent commuting integrals.  Solutions are generated as logarithmic
derivatives of so called $\tau$-functions.  
The partition functions (\ref{partfun}) 
are these  $\tau$-functions.

Using this fact we showed that the
leading order (as $N\to \infty$)
 of a second logarithmic derivative of the partition
function was the generating function of the higher Catalan numbers.  
\begin{thm}[Ercolani-McLaughlin-Pierce] \label{theorem1.2}
\begin{equation*}
z(-c_\nu t) + \mathcal{O}\left(\frac{1}{N^2}\right) 
= \frac{1}{2N^2} \frac{\partial^2}{\partial t_1^2}
\log\left[ \frac{1}{Z_N(0, 0)} Z_N( N^{1/2} t_1, N^{1-\nu} t )
  \right]_{t_1 = 0} \, .
\end{equation*}
\end{thm}

The paper \cite{EMP} presents two independent
proofs of Theorem \ref{theorem_1} in the case when $g=0$: one is
a direct calculation of the 
generalization of the method of steepest descent for computing the
leading order of asymptotics of matrix integrals, 
the other uses the fact that
the partition function (\ref{partfun})
 is a $\tau$-function of the Toda lattice
hierarchy.  
Theorem \ref{theorem1.2}  is found in  either proof of Theorem
\ref{theorem_1} in the $g=0$ case.
The general result from this paper is that $e_g(t)$ is explicitly
computable as a function of $z(s)$.  
Therefore the higher Catalan
numbers play a central role in the combinatorics of maps.

\subsection{Results}

We will summarize our results for enumerations of genus 0 and genus 1
maps  here.  The generating
function of the number of genus 0 maps is 
\begin{equation*}
e_0(-t_{2\nu}) = \sum_{j=0}^\infty \kappa^{(2\nu)}_0(j)
\frac{t_{2\nu}^j}{j!} ,
\end{equation*}
while the generating function for genus 1 maps is 
\begin{equation*}
e_1(-t_{2\nu}) = \sum_{j=0}^\infty \kappa^{(2\nu)}_1(j)
\frac{t_{2\nu}^j}{j!} .
\end{equation*}

We will give a combinatoric proof of the theorem from \cite{EMP}:
\begin{thm}[Ercolani-McLaughlin-Pierce] \label{thm1}
The Taylor coefficients of $e_0$ are 
\begin{equation*}
\kappa_0^{(2\nu)} = \frac{c_\nu^j (\nu j - 1)! }{( (\nu-1) j + 2)!} ,
\end{equation*}
and those of $e_1$ are
\begin{equation*}
\kappa_1^{(2\nu)} = \frac{ (j-1)! c_{2\nu}^j }{12} \sum_{k=1}^j
(\nu-1)^k \binom{\nu j}{j-k} .
\end{equation*}
\end{thm}

In section \ref{calc} we will review the combinatoric arguments
showing that the higher Catalan numbers are generated by $z(s)$.  Both
the path counting problems analyzed in section \ref{calc} and the
queuing problem mentioned are from Yaglom and Yaglom \cite{YY}, and
Sato and Cong \cite{Sa}.  The path counting problems have been
cosmetically changed, but the method of solution is still the same. 
In section \ref{e0_sec} we give the origin of theorems like
\ref{theorem_1} and our general result.   
In section \ref{e0_sec2} we give a combinatoric calculation starting
from Theorem \ref{theorem_1} which culminates in a formula for the
number of planar maps with $j$ vertices of degree $2\nu$.  In section
\ref{e1_sec} we carry out a similar calculation for counting genus 1
maps, giving a proof of Theorem \ref{thm1}.

\section{\label{calc} Calculation of the Higher Catalan
  Numbers} 

The higher Catalan numbers satisfy a recursion relation which will tie
together four types of combinatoric problems:  the Taylor coefficient
problem (\ref{z(s)_taylor}), the polygon division problem, and two
problems we are about to introduce: a path counting problem and a
queuing problem.  The path counting problem is directly solvable by
enumeration.  This is then connected to the queuing problem and we
show that both
 satisfy the same recursion (\ref{recur})
   as $\zeta_j^{(\nu)}$ and the polygon division problem.

We give here a modification of the calculation outlined in Yaglom and Yaglom
\cite{YY}, and Sato and Cong \cite{Sa}.  
In this section we will show that the higher Catalan numbers, 
\begin{equation*}
\zeta_j^{(\nu)} = \frac{1}{j} \binom{\nu j}{j-1},
\end{equation*}
count the number of right 1, left $(\nu-1)$ walks  on the 1-Dimensional
lattice which go from 0 to 0 in $\nu j$ steps which stay to the right of 0.  

Note that this 1-Dimensional counting problem is identical to
counting the number of $(1,1)$ and $(1, -(\nu-1))$ paths on the 2-Dimensional
Integer Lattice from $A=(0, 0)$ to $B=(\nu j, 0)$ which stay above
$y=0$, we call these higher Dyck paths.

As an example of how to count such paths we have the lemma
\begin{lem}
The number of $(1,1)$ and $(1, -(\nu-1))$ paths from $P= (p_1, p_2)$ to $Q
=(q_1, q_2)$ ( assume $ p_1 < q_1$) on the 2-Dimensional integer lattice is 
\begin{equation*}
N_{PQ} = \binom{q_1-p_1}{d}\,,
\end{equation*}
where $n, u,$ and $d$ solve the system of equations
\begin{align*}
p_1 + u + d = q_1 \,, \\
p_2 + u - (\nu-1) d = q_2.
\end{align*}
If no integer solutions exist then there is no such path connecting these
points.  
\end{lem}

To count the number of $(1,1)$, $(1, -(\nu-1))$ paths from $A$ to $B$ which
stay above the $x$-axis we will count the number of such paths which pass below
the $x$-axis.  A path of this type going below the $x$-axis will go through
one of the points 
\begin{equation*}
D_k = \left( \nu k - 1, -1 \right).
\end{equation*}
Let $A'=(0, -\nu)$.  
\begin{lem}
The number of $(1,1)$, $(1, -(\nu-1))$ paths from $A$ to $B$ which pass below
the $x$-axis is 
\begin{equation*}
(\nu-1) N_{A'B} = (\nu-1)\binom{\nu j}{j-1}.
\end{equation*}
\end{lem}

Begin the proof of this lemma by computing
\begin{equation*}
N_{AD_k} = \binom{\nu k -1}{k} 
\end{equation*}
and
\begin{equation*}
N_{A'D_k} = \binom{\nu k -1}{k-1} = \frac{1}{\nu-1} \binom{\nu k -1}{k} =
\frac{1}{\nu-1} N_{AD_k}.
\end{equation*}
Now consider the number of paths from $A$ to $B$ which go through $D_1$
\begin{equation*}
N_{AD_1} N_{D_1 B} = (\nu-1) N_{A' D_1} N_{D_1 B},
\end{equation*}
which is $(\nu-1)$ times the number of paths from $A'$ to $B$ through $D_1$.  
Then inductively assume that the number of paths from $A$ to $B$ through $D_l$
is $(\nu-1)$ times the number of paths from $A'$ to $B$ through $D_l$ for
$l<k$.  
Consider the number of paths from $A$ to $B$ which go through $D_k$ but not
$D_l$ for $l<k$,
\begin{equation*}
N_{AD_k} N_{D_k B} - N_{A D_{k-1}} N_{D_{k-1} D_k} N_{D_k B} - \dots 
= (\nu-1) N_{A' D_k} N_{D_k B} - (\nu-1) N_{A' D_{k-1}} N_{D_{k-1} D_k} N_{D_k
  B} - \dots.
\end{equation*}
\qed

The number of $(1,1)$, $(1, -(\nu-1))$ paths from $A$ to $B$ which stay above
the $x$-axis is 
\begin{equation} \label{star}
\binom{\nu j}{j} - (\nu-1) \binom{\nu j}{j-1} = \frac{1}{j} \binom{\nu
  j}{j-1}, 
\end{equation}
which is the $j$'th $\nu$-higher Catalan number.   

Consider the queuing problem:
How many ways are there for $(\nu-1) j$ customers with $1$-dollar bills
and $j$ customers with $\nu$-dollar bills to line up to buy $1$-dollar widgets
so that each customer receives correct change.  
In the path counting problem, in order for a step to the left to be 
placed on the path and stay to the right of 0 there must be
$(\nu-1)$ corresponding right steps that have come before on the
path.  
In the queuing problem, in order for a customer with a $\nu$-dollar bill to
receive correct change there must be $(\nu-1)$ corresponding customers with
$1$-dollar bills in line ahead of the $\nu$-dollar bill.

To see that the numbers (\ref{star}) agree with the coefficients
(\ref{High_Cat}) of
$z(t)$ we argue that the number of right 1, left $(\nu-1)$ paths  which stay
to the right of 0 satisfy the recursion relation (\ref{recur}).  As argued
above, to each left step in the path there must correspond $(\nu-1)$ steps to
the right.  The path begins with a step to the right.  This step corresponds,
together with $(\nu-2)$ other right steps, to a left step later on.  This
collection of $\nu$ steps divides the entire $\nu j$ steps into $\nu$
sub-paths of the type: steps right by 1, left by $(\nu-1)$, from $k$ to $k$
which stay to the right of $k$.  Therefore the recursion relation
(\ref{recur}) is satisfied.

Now we argue
 that the number of ways to divide a marked $((\nu-1)j +2)$-sided polygon
into $(\nu+1)$-sided polygons using edges and diagonals satisfies the
recursion relation (\ref{recur}).  Take the sub-polygon which has the marked
outside edge as an edge.  This sub-polygon divides the $((\nu-1)j +2)$-sided
polygon into $\nu$ polygons, therefore the recursion relation
 (\ref{recur}) is satisfied by this problem as well.  
Therefore these numbers are also the higher Catalan numbers.

\section{The Enumeration of Planar Maps \label{e0_sec}}

Our interest in these combinatoric problems began with the function
$z(s)$ defined by (\ref{z(s)_taylor}).  This function appeared as an auxiliary
function in explicit calculations of the asymptotic expansion of the
partition function of large random matrices.  It is interesting that a
generating function as classical as $z(s)$ appeared in this setting as
the terms of this asymptotic expansion are themselves generating
functions for a combinatoric problem.  The Taylor coefficients of the
terms in the asymptotic expansion enumerate maps on genus $g$
surfaces.  

In this section we highlight the information contained in the asymptotic
expansion of the partition function (\ref{partfun})  for large matrices.  
 We have detailed data about the higher genus
problems \cite{EMP}, however here we will concentrate on the genus
zero (or planar) setting.

Recall that the partition function of random matrices we consider is
(\ref{partfun}): 
\begin{equation*}
Z_N(t) = \int_{\mathcal{H}_N} \exp\left( -N \mbox{Tr}\left( \frac{1}{2} M^2 +
t M^{2\nu} \right) \right) dM,
\end{equation*}
where now we take $t_1 = 0$.  
  Ercolani and
McLaughlin \cite{EM} showed that 
$\log\left( Z_N(t)/Z_N(0) \right)$
possesses an asymptotic expansion as $N \to \infty$ 
inside a non-trivial $t$ domain, 
\begin{equation} \label{asympt}
\frac{1}{N^2} \log\left(\frac{Z_N(t)}{Z_N(0)}\right) = 
e_0(t) + \frac{1}{N^2} e_1(t) + \frac{1}{N^4} e_2(t) + \dots
\end{equation}
where $e_g(t)$ is analytic in a neighborhood of $t=0$ 
and is a counting function for
genus $g$ maps with $2\nu$-degree 
vertices  The idea of counting maps in this way originates with random
matrix models of quantum field theories \cite{BIZ, FGZ, Wi}.

In this paper we will concentrate our attention on the genus $0$ and
genus $1$ terms in (\ref{asympt}).  The genus $0$ term is 
\begin{equation*}
e_0(t) = \sum_{j=1}^\infty \kappa_0^{(2\nu)}(j) \frac{(-t)^j}{j!},
\end{equation*}
where the 
$\kappa_0^{(2\nu)}(j)$ are the number of genus $0$ maps with $j$ vertices
of degree $2\nu$.  
These  numbers are calculated in 
\cite{EMP} by explicit computation of the contour integral
\begin{equation} \label{contour}
\kappa_0^{(2\nu)}(j) = \frac{j!}{2\pi i} \oint t^{-j-1} e_0(-t) dt\,,
\end{equation}
where the contour encircles $t=0$ and $e_0$ is given in Theorem
\ref{theorem_1}.  To evaluate (\ref{contour}) one rewrites it as
an integral with respect to $z = z(c_\nu t)$.  
The result is the first part of Theorem \ref{thm1}.

\section{\label{e0_sec2} Computation of the Taylor Coefficients of $e_0$}

We will now compute the Taylor coefficients
of $e_0$ from combinatoric arguments only.  
In section \ref{calc} we used combinatoric techniques to
compute the Taylor coefficients of $z(s)$.  Theorem \ref{theorem_1} gives
$e_0(-t)$ as a function of $z(s)$.  The calculation of the Taylor coefficients
of $e_0(-t)$ can be done from this theorem utilizing contour integration and
a clever change of variables.  Our goal in this note 
is to complete this calculation using
only combinatoric arguments.  To that end, we need to derive formulas for the
 coefficients of $\log(z(s))$ and $(z(s)-1)^2$.

\subsection{Taylor Coefficients of powers of $(z(t)-1)$ and $\log(z(t))$}

The coefficients of $(z(s)-1)^i$ are given by 
\begin{equation*} %\label{zsqr}
\eta_j^{(\nu, i)} = \sum_{j_1 + j_2 \dots + j_i = j \\
j_n > 0 } \zeta_{j_1}^{(\nu)}
\zeta_{j_2}^{(\nu)} \dots \zeta_{j_i}^{(\nu)}, \; j\geq i .
\end{equation*}
The $\eta_j^{(\nu, i)}$ 
count the number of right 1, left $(\nu-1)$ paths from 0 to 0 in
$\nu j$ steps which stay to the right of 0 and return to 0 at least $i-1$
times in
between the ends of the path.
This is equivalent to the queuing problem:  how many ways can $j$ customers
with $\nu$-dollar bills and $(\nu-1) j$ customers with $1$-dollar bills form
$i$ lines to buy 1 dollar widgets so that each customer receives exact
change.  

This multi-line 
queuing problem is equivalent to:  how many ways can $j-i$ customers with
$\nu$-dollar bills and $(\nu-1)j + 1$ customers with 1-dollar bills form one
line to buy 1 dollar widgets so that each customer receives exact change.    
This is done by shifting the lines together.  

First, we notice that at the back of the first line is a $\nu$-dollar
bill.  We shift the problem in the following way:  we change the $\nu$-dollar
bill at the back of first line into a 1-dollar bill, then adjoin the second
line to the back of the first.  
The first line now has $j-1$ customers
with $\nu$-dollar bills and $(\nu-1) j + 1$ customers with 1-dollar bills (a
$\nu$-dollar bill has been changed into a 1-dollar bill).  
We are now in the case of $i-1$ lines and may shift again.
After $i-1$ repetitions of this process we have a single line 
arranged in such a way that the last customer has a
$\nu$-dollar bill, we remove this customer from the process.

To show that the two counting problems are equivalent we show that the
one-line queuing problem above can be transformed back into the multi-line
problem.  We begin with a line of $j-i$ customers with $\nu$-dollar bills and
$(\nu-1)j + 1$ customers with $1$-dollar bills in a single line to buy 1
dollar widgets so that each customer will receive exact change.  First add a
$\nu$-dollar bill to the back of the line.  Starting from the back of the line
and moving forward we will form sets of blocks of $\nu$ customers.  Each block
has a customer with a $\nu$-dollar bill and $\nu-1$ customers with 1-dollar
bills.  When we reach a customer with a 1-dollar bill for which there is no
corresponding $\nu$-dollar bill behind, we have completed a line, and
we remove it from the first line.  
We then turn the 1-dollar bill at the back of the
first line into a $\nu$-dollar bill and proceed as before.  Thus transforming
us to the multi-line queuing problem.

The path counting problem which corresponds to this shifted queuing problem
is:  how many right 1, left $(\nu-1)$ paths, from 0 to $(\nu-1)i + 1$,  
in $\nu j-i+1$ steps, are there which stay to the right of 0.
This is equivalent to the problem:  count the number of
 $(1, 1)$, $(1, -(\nu-1))$ paths
from $(0, 0)$ to $(\nu j - i +1, (\nu-1)i + 1)$, which stay above zero.  
These paths are counted in the same way as the ones in section 2:
\begin{equation}\label{zsqr_form}
\eta_j^{(\nu, i)}  = \frac{i}{j} \binom{\nu j}{j-i}, \; j\geq i .
\end{equation}

The coefficients of $\log(z(s))$ are computed directly using the Taylor
expansion of the logarithm function:
\begin{equation*}
\log\left( z_0(s) \right) = \sum_{j=1}^\infty l^{(\nu)}_j s^j =
\sum_{i=1}^{\infty} (-1)^{i} 
\frac{\left( z_0(s) -  1\right)^i}{i} .
\end{equation*}
We find that the coefficients are
\begin{equation} \label{logz_form}
l_j^{(\nu)} = \sum_{i=1}^j (-1)^i \frac{ \eta_j^{(\nu, i)}}{i} = 
\sum_{i=1}^j (-1)^i \frac{1}{j} \binom{\nu j}{j-i} =
\frac{1}{j} \binom{\nu j-1}{j-1}\,.
\end{equation}

\subsection{Enumeration of planar maps}

Combine the results (\ref{zsqr_form}) and (\ref{logz_form}) with Theorem
\ref{theorem_1}:
\begin{equation*}
e_0(-t) = \mu (z(c_\nu t) - 1) ( z(c_\nu t) - r) + \frac{1}{2} \log( z(c_\nu
t) ) ,
\end{equation*}
and we find that
\begin{align*}
\kappa_0^{(2\nu)}(j) &= j! c_\nu^j \left[ - (r-1) \mu \zeta_j^{(\nu)} + \mu
  \eta_j^{(\nu, 2)} + \frac{1}{2} l_j^{(\nu)} \right] 
\\
&= c_\nu^j \frac{( \nu j - 1)!}{ ( (\nu-1) j + 2)! } \,.
\end{align*}
This proves the first part of Theorem \ref{thm1}.

\subsection{Additional enumerative results for functions of $z(s)$}

Formulas similar and equivalent to 
(\ref{zsqr_form}) have been derived for other
functions of $z(s)$ \cite{Go56, Go76}.  

\begin{thm}[Polya-Szeg\"{o}-Gould \label{psg} ]
Let $z(s)$ be the generating function for the Higher Catalan Numbers.  
Then 
\begin{equation*}
z^\alpha = \sum_{j=0}^\infty A_j(\alpha, \nu) s^j,
\end{equation*}
with 
\begin{equation}\label{Aj}
A_j(\alpha, \nu) = \frac{\alpha}{(\alpha+\nu j)} \binom{ \alpha + \nu
  j}{j}.
\end{equation}
In addition
\begin{equation*}
\frac{z^{\alpha+1}}{\nu - (\nu-1) z} = \sum_{j=0}^\infty
\binom{\alpha+\nu j}{j} s^j.
\end{equation*}
\end{thm}

One can check that formula (\ref{Aj}) is equivalent to
(\ref{zsqr_form}).

\section{\label{e1_sec} The enumeration of genus one maps}

In the paper \cite{EMP} the authors show that there is a general
construction of the functions $e_g(t_{2\nu})$ in the asymptotic
expansion (\ref{asympt}) as functions of $z(c_{2\nu} t_{2\nu})$.  
For example, 
we find that in the case of
genus one maps
\begin{thm}[Ercolani-McLaughlin-Pierce \label{theorem2}]
In the neighborhood of analyticity, the second term in the asymptotic
expansion (\ref{asympt}) is given by 
\begin{equation} \label{e1form}
e_1(-t) = -\frac{1}{12} \log\left( \nu - (\nu-1) z(c_{2\nu} t) \right) ,
\end{equation}
where we choose the principal branch of the logarithm.  
\end{thm}

Our goal is to derive a formula for the Taylor coefficients of
$e_1(-t_{2\nu})$ for general $\nu$.  This is a straightforward
calculation: expand formula (\ref{e1form}) as a powers series in
$(z( c_{2\nu} t) - 1)$;
\begin{equation}\label{e1step1}
e_1(-t) = \sum_{k=1}^\infty \frac{1}{12} \frac{(\nu-1)^k}{k} \left(
z(c_{2\nu} t) - 1 \right)^k .
\end{equation} 
Then insert the power series of $(z - 1)^k$ into (\ref{e1step1}):
\begin{equation}\label{e1step2}
e_1(-t) = \sum_{k=1}^\infty \sum_{j=k}^\infty \frac{1}{12}
\frac{(\nu-1)^k}{k} \eta^{(\nu, k)}_j t^j = 
\sum_{k=1}^\infty \sum_{j=k}^\infty \frac{(\nu-1)^k}{12 j} \binom{ \nu
  j}{j-k} t^j ,
\end{equation}
where the $\eta_j^{(\nu, k)}$ were computed in (\ref{zsqr_form}).
Switching the order of summation in (\ref{e1step2}) we find that
\begin{equation}\label{unstuck}
\kappa_1^{(2\nu)}(j) = \frac{(j-1)! c_{2\nu}^j}{12} \sum_{k=1}^j
(\nu-1)^k \binom{\nu j}{j-k} .
\end{equation}
This proves the second part of Theorem \ref{thm1}.  

One is tempted to try to use the last part of Theorem \ref{psg} to
find a more concise formula for $\kappa_1^{(2\nu)}(j)$.  However using
the relation 
\begin{equation*}
z' = \frac{ z^{\nu+1}}{\nu - (\nu-1) z}, 
\end{equation*}
derived from the implicit definition of $z$ by formula (\ref{constraint}),
 we see that the derivative of $e_1(-t)$ with respect to $t$ is 
\begin{equation}\label{stuck}
e_1'(-t) = \frac{(\nu-1)}{12} \frac{ z(c_{2\nu} t)^{\nu+1} }{(\nu -
  (\nu-1) z(c_\nu t) )^2}. 
\end{equation}
There does not seem to be a more concise formula for the coefficients
of (\ref{stuck}) than that found directly from $j$ times
$\kappa_1^{(2\nu)}$ in (\ref{unstuck}).

\section{Conclusion}

The higher Catalan numbers appeared naturally in the study of the
enumeration of maps embedded on Riemann surfaces.  We have gathered
here a number of interesting facts about these numbers and the
combinatoric problems related to them.    
It would be interesting to have a combinatoric argument which gives
the form of Theorems \ref{theorem_1} and \ref{theorem2}.  

Further work will study the fine structure of the
asymptotic expansion of the partition function $Z_N$ when multiple
time evolutions are involved. Additionally other partition functions
over different families of Random Matrices encode
similar combinatoric data.  
For example, in the case when 
\begin{equation*}
Z_N = \int_{\mathcal{S}_N} \exp\left[ -N \mbox{Tr}\left( \frac{1}{4}
  M^2 + V(M) \right) \right] dM,
\end{equation*}
 where the integral is taken over the space of $N \times N$ symmetric
 matrices,  we find that the terms in the asymptotic expansion of
 $\log\left(Z_N\right)$ give generating functions for the number of
 unoriented  maps partitioned by the Euler characteristic associated with the
 embedding of the map.  

What connections these problems will have to the higher Catalan
numbers or to other interesting (and classical) combinatoric problems
remains to be seen.

\bibliographystyle{amsplain}

\begin{thebibliography}{9}
\bibitem{Ba}
D. F. Bailey, \textit{Counting Arrangements of 1's and -1's}, Mathematics
Magazine, \textbf{69},  (1996)  128-131.
\bibitem{BIZ}
D. Bessis, C. Itzykson, and J.B. Zuber, \textit{Quantum field theory
  techniques in graphical enumeration}, Advances in Applied Mathematics,
\textbf{1}, (1980), no. 2, 109-157.
\bibitem{BB}
J. Borwein and D. Bailey, \textit{Mathematics by Experiment: Plausible
  Reasoning in the 21st Century}  Natick, MA:  A. K. Peters, 2003.  
\bibitem{Br}
R.A. Brualdi, \textit{Introductory Combinatorics, 4th ed}, New York:  Elsevier
(1997). 
\bibitem{CG}
J.H. Conway and R.K. Guy, \textit{The book of numbers}, New York:
Springer-Verlag, (1996)  96-106.
\bibitem{FGZ}
P. DiFrancesco, P. Ginsparg, and J. Zinn-Justin, \textit{2D gravity and random
  matrices}, Physics Reports, \textbf{254} , (1995), 1-133.
\bibitem{Do}
H. D\"{o}rrie, \textit{Euler's problem of Polygon Division}, \S7 in
\textit{100 Great Problems of Elementary Mathematics: Their History and
  Solutions}, New York: Dover (1965), 21-27.
\bibitem{Go56}
H.W. Gould, \textit{Some generalizations of Vandermonde's
  convolution}, American Math. Monthly, \textbf{63}, (1956) 84-91.
\bibitem{Go76}
H.W. Gould, \textit{Research bibliography of two special number
  sequences}, Mathematicae Monongaliae, Dept. of Math., West Virginia
Univ.  (1971, revised 1976).
\bibitem{EM}
N. M. Ercolani and K. D. T.-R. McLaughlin, \textit{Asymptotics of the
  partition function for random matrices via Riemann-Hilbert techniques and
  application to graphical enumeration}, International mathematics research
notes \textbf{14} (2003), 755-820.
\bibitem{EMP}
N. M. Ercolani, K. D. T.-R. McLaughlin, and V. Pierce, 
\textit{Random
  matrices, graphical enumeration and the continuum limit of the Toda
  lattices}, math-ph/0606010 (2006).
\bibitem{Ho}
R. Honsberger, \textit{Mathematical Gems I}, Washington, DC:  Mathematical
Association of America,  (1973), 130-134.
\bibitem{Kl}
D.A. Klarner, \textit{Correspondences between plane trees and binary
  sequences}, Journal of Combinatoric Theory, \textbf{9} (1970)
401-411.
\bibitem{LaYe90}
J. Labelle and Y.-N. Yeh, \textit{Generalized Dyck paths}, 
Discrete Mathematics, \textbf{82} (1990) 1-6.
\bibitem{MW}
M.E. Mays and J. Wojciechowski, \textit{A Determinant Property of Catalan Numbers},
Discrete mathematics, \textbf{211}, (2000) 125-133.  
\bibitem{Po}
G. P\'{o}ya, \textit{On Picture-Writing}, American Mathematical Monthly
\textbf{63}, (1956)  689-697.
\bibitem{PS}
G. Polya and G. Szeg\"{o}, \textit{Aufgaben und Lehrs\"{a}tze aus der
  Analysis}, Vol I, Springer-Verlag, Berlin (1925).
\bibitem{Sa}
M. Sato, and T. T. Cong, \textit{The number of minimal lattice paths
  restricted by two parallel lines.}  Discrete mathematics,
\textbf{43}, (1983), 249-261.
\bibitem{Wi}
E. Witten, \textit{Two-dimensional gravity and intersection theory on moduli
  space.}  Surveys in Differential Geometry  \textbf{1},  (1991), 243-310.
\bibitem{YY}
A. M. Yaglom and I. M. Yaglom, \textit{Challenging mathematical problems with
  elementary solutions}, vol. I:  Combinatorial Analysis and Probability
Theory, Holden-Day, 1964.
\bibitem{Zv}
A. Zvonkin, \textit{Matrix integrals and map enumeration: an
  accessible introduction.}  Mathematical and computational modeling,
\textbf{26} (1997) 281-304.
\end{thebibliography}

\end{document}